\documentclass[a4paper,11pt]{article}
	\usepackage{amsmath,amsthm,amssymb,amsfonts,epsfig,epstopdf,titling,url,array,enumerate,cite,mathrsfs,upgreek,authblk,scrextend, graphicx, float}
	\usepackage[bindingoffset=0.2in,left=1in,right=1in,top=1in,bottom=1.5in,footskip=.25in]{geometry}
	\theoremstyle{plain}
	\newtheorem{thm}{Theorem}[section]
	\providecommand{\keywords}[1]{\begin{addmargin}[28pt]{28pt}\noindent\textbf{Keywords:} #1 \end{addmargin}}
	\newtheorem{lma}[thm]{Lemma}

	\theoremstyle{definition}
	\newtheorem{dfn}[thm]{Definition}
	\newtheorem{eg}[thm]{Example}
	\newtheorem{rem}[thm]{\textit{Remark}}
	\providecommand{\ams}[1]{\begin{addmargin}[28pt]{28pt}\noindent\textbf{Mathematics Subject Classification:} #1\end{addmargin}}
	
	\title{Some basic results on fuzzy strong $\phi$-b-normed linear spaces }
\author[1]{Abhishikta Das}
\author[2,*]{T. Bag}
\author[3]{S. Chatterjee}
\affil[1,2,3]{Department of Mathematics, Siksha-Bhavana,	\authorcr	Visva-Bharati, Santiniketan-731235, Birbhum, West-Bengal, India  
\authorcr 	E-mail\textsuperscript{1}: abhishikta.math@gmail.com	\authorcr E-mail\textsuperscript{2,*}: tarapadavb@gmail.com 
\authorcr  E-mail\textsuperscript{3}: shayani.mathvb10@gmail.com 	}

	\date{}

\begin{document}
		\maketitle
	\begin{abstract}
	
		\noindent
In this paper, definition of fuzzy strong $\phi$-b-normed linear space is given. Here the scalar function $|c| $ is replaced by a general function 
$ \phi(c) $ where $ \phi $ satisfies some properties. Some basic results on finite dimensional  fuzzy strong $\phi$-b-normed linear space are studied.
	\end{abstract}
	\ams{54A40, 03E72}
	\section{Introduction } 
	The concept of a fuzzy set was introduced initially by Zadeh\cite{[14]} in 1965. Since then, many authors have expansively developed the
	 theory of fuzzy sets. Osmo Kaleva\cite{[17]}, Kramosil and Michalek\cite{[16]}, Georage and Veeramani\cite{[15]} et al. introduced the
	  concept of fuzzy metric spaces in different approaches. on the other hand, Katsaras\cite{[19]}, Felbin\cite{[5]}, Cheng and Mordeson\cite{[4]},
	  Bag and Samanta\cite{[1]} gave the definition of fuzzy normed linear spaces in different way. \\
	  Recently different types of generalized metric  as well as norm (viz. 2-metric\cite{[26]},   b-metric\cite{[9]}, strong b-metric\cite{[10]},
	   G-metric\cite{[23]}, 2-norm\cite{[27]}, G-norm\cite{[28]}, etc.) and consequently generalized fuzzy metric and fuzzy norm
	   (viz. fuzzy b-metric\cite{[24]}, strong fuzzy b-metric\cite{[25]}, fuzzy cone metric\cite{[21]},
	    fuzzy cone norm\cite{[22]},  G-fuzzy norm\cite{[29]},  etc. ) have been introduced in different approaches. \\
	    In 2018, Oner\cite{[25]} introduced the concept of fuzzy strong b-metric spaces and developed some topological results in such spaces. Following 
	    this definition of fuzzy strong b-metric spaces, in this paper we give a definition of fuzzy strong $\phi$-b-normed linear space whose induced
	    fuzzy metric is Oner type.	In fuzzy normed linear space, scalar multiplication is given by $ N(cx,t) = N(x, \frac{t}{|c|}) $. But in our
	     definition of  fuzzy strong $\phi$-b-norm, scalar multiplication is given by $ N(cx,t) = N(x, \frac{t}{\phi(c)}) $ where $ \phi $ 
	     is a real valued function satisfying some properties. We study some results on finite dimensional fuzzy strong $ \phi$-b-normed linear spaces.\\
	     The organization of the paper is in the following.\\
	      Section 2 consists some preliminary results. In Section 3, we introduce a definition of fuzzy strong $ \phi$-b-norm by using 
	      a special function $ \phi $ in general t-norm settings and illustrate by examples.  Some basic results of finite dimensional fuzzy 
	      strong $ \phi$-b-normed linear spaces  are established in Section 4. 
\section{ Preliminaries}
In this section some definitions and  results are collected which are used in this paper.
	  \begin{dfn} \cite{[3]}
 A binary operation $ * : [0 , 1] \times [0 , 1] \rightarrow [0 , 1] $ is called a  $t$-norm if it satisfies the following conditions:
 \begin{enumerate}[(i)]
		\item  $ * $ is associative and commutative; 
		\item $ \alpha ~ * ~ 1 = \alpha, ~ \forall \alpha \in [0 , 1] $;  
	\item $ \alpha * \gamma \leq  \beta * \delta ~ $  whenever $ \alpha \leq  \beta  $ and $  \gamma \leq   \delta,  ~\forall \alpha, \beta,
	\gamma, \delta \in [0 , 1] $.
	\end{enumerate}
If $ * $ is continuous then it is called continuous $t$-norm.
\end{dfn}
 The following are examples of some t-norms. 
\begin{enumerate}[(i)]
		\item  Standard intersection: $ \alpha *   \beta = \min \{ \alpha ,   \beta \} $.
		\item   Algebraic product: $ \alpha *   \beta = \alpha    \beta $.
	\item Bounded difference: $ \alpha *   \beta = \max \{ 0,  \alpha +  \beta -1 \} $.
	\end{enumerate}
\begin{dfn} \cite{[2]}
 	Let $ X $ be a linear space over a field $ \mathbb{F} $. A  fuzzy subset $ N $ of $ X \times \mathbb{R}$ 
 	 is called fuzzy norm on $  X $ if  $ \forall x, y \in X $ the following conditions hold:
 	\begin{enumerate}[(N1)] 
 		\item $ \forall t \in {\mathbb{R}}  $ with $ t \leq 0, ~ N(x,t)=0 $; 
 		\item $ (\forall t \in {\mathbb{R}}, ~ t>0, ~ N(x,t) = 1)   \iff  x = \theta $;
 		\item $ \forall t \in {\mathbb{R}}, $ and  $ c \in {\mathbb{R}}, ~ t>0, ~ N(cx,t) = N(x, \frac{t}{|c|}) $;
 		\item $ \forall s, t \in {\mathbb{R}}, ~ N(x+y, s+ t) \geq N(x,s) * N(y,t) $;
 		\item $ N(x, \cdot) $ is a non-decreasing function of $ t $ and $ \underset{t \rightarrow \infty }{lim} N(x,t) = 1 $. 
 	\end{enumerate}
 	Then the pair  $ ( X, N ) $ is called  fuzzy normed linear space. 	
  \end{dfn}
\begin{dfn} \cite{[1]}
Let $ ( X, N ) $  be a fuzzy normed linear space.
\begin{enumerate}[(i)]
\item A sequence $ \{ x_n \} $ is said to be convergent if $ \exists x \in X $ 
such that $ \underset{n \rightarrow \infty }{lim} N(x_n - x, t) = 1, ~ \forall t > 0 $. Then  $ x $ is called the limit of the sequence $ \{ x_n \} $
and denoted by $ lim x_n $.
\item A sequence $ \{ x_n \} $ in a fuzzy normed linear space $ ( X, N ) $ is said to be Cauchy if \\
 $ \underset{n \rightarrow \infty }{lim} N(x_{n+p} - x_n, t) = 1, ~ \forall t > 0 $ and $ p = 1,2 \cdots $.
 \item  $ A \subseteq  X $ is said to be closed if for any sequence $ \{ x_n \} $ in $ A $ converges to $ x \in A $. 
 \item $ A \subseteq  X $ is said to be the closure of $ A $, denoted by $ \bar{A} $ if for any $x \in \bar{A} $, there is a sequence
 $ \{ x_n \} \subseteq A $ such that $ \{ x_n \} $ converges to  $ x $. 
 \item $ A \subseteq  X $ is said to be compact if any sequence $ \{ x_n \} \subseteq A $  has a subsequence converging to
an element of $ A $.
 \end{enumerate}
\end{dfn}
\begin{dfn} \cite{[20]}
Let $ ( X , N) $ be a fuzzy normed linear space. 
\begin{enumerate} [(i)]
\item A set $ B(x, \alpha, t), ~ 0 < \alpha  < 1 $ is defined  as $ B(x, \alpha, t) = \{ y : N(x - y , t) > 1 - \alpha \} $.
\item $ \tau = \{ G \subseteq X : x \in G, ~ \exists t > 0 ~ \text{and} ~ 0 < \alpha < 1 ~ \text{such that} ~  B(x, \alpha, t) \subset G \} $
is a topology on $ (X, N ) $.
\item Members of $ \tau $ are called open sets in $ (X, N ) $.
\end{enumerate}
\end{dfn}
\begin{dfn} \cite{[2]}
 A subset $ B $ of a fuzzy normed linear space $ ( X, N ) $  is said to be fuzzy bounded if for each $ r, ~ 0 < r < 1, ~\exists t > 0 $
  such that $ N(x, t) > 1, ~ \forall x \in B $.
\end{dfn} 
	\begin{lma}\cite{[2]}
	Let $(X , N) $ be a fuzzy normed linear space and $ N(x , \cdot)  (x \neq 0)$ .
If the set $ M = \{ x : N(x , 1) > 0 \} $ is compact then $X$ is finite dimensional.
\end{lma}
\section{Fuzzy strong $\phi$-b-normed linear space}
%
In this section we give the definition of fuzzy normed linear space in  a new approach.
 \begin{dfn}
	  Let  $ \phi $ be a function defined on  $ {\mathbb{R}} $ to $ {\mathbb{R}} $ with the following properties 
	\begin{enumerate}[($\phi$1)]
		\item $ \phi (-t) = \phi (t), ~ \forall t \in {\mathbb{R}}  $; 
		\item $ \phi (1) = 1 $;
		\item $ \phi $ is strictly increasing and continuous on $ (0, \infty) $;
		\item $ \underset{\alpha \rightarrow 0 }{lim} \phi(\alpha) = 0 $ and $ \underset{\alpha \rightarrow \infty }{lim} \phi(\alpha) = \infty $.
	\end{enumerate}
	\end{dfn}
The followings are examples of such functions.
 \begin{enumerate}[(i)]
		\item $ \phi(\alpha) = | \alpha |,  ~ \forall \alpha \in {\mathbb{R}}  $
		\item $ \phi(\alpha) = | \alpha |^p, ~\forall \alpha \in {\mathbb{R}}, ~p \in \mathbb{R}^+  $
		\item $ \phi(\alpha) = \frac{2 \alpha^{2n}}{|\alpha| + 1}, ~ ~\forall \alpha \in {\mathbb{R}}, ~ n \in \mathbb{N} $
	\end{enumerate}
 \begin{dfn}\label{dfn1}
 	Let $X$ be a linear space over a field $ \mathbb{F} $ and $ K \geq 1 $  be a given real number. A  fuzzy subset $ N $ of $ X \times \mathbb{R}$ 
 	 is called fuzzy strong $\phi$-b-norm on $  X $ if  $ \forall x, y \in X $ the following conditions hold:
 	\begin{enumerate}[(bN1)] 
 		\item $ \forall t \in {\mathbb{R}}  $ with $ t \leq 0, ~ N(x,t)=0 $; 
 		\item $ (\forall t \in {\mathbb{R}}, ~ t>0, ~ N(x,t) = 1)  $ iff $ x = \theta $;
 		\item $ \forall t \in {\mathbb{R}}, ~ t>0, ~ N(cx,t) = N(x, \frac{t}{\phi(c)} ) $ if $\phi(c) \neq 0 $;
 		\item $ \forall s, t \in {\mathbb{R}}, ~ N(x+y, s+ Kt) \geq N(x,s) * N(y,t) $;
 		\item $ N(x, \cdot) $ is a non-decreasing function of $ t $ and $ \underset{t \rightarrow \infty }{lim} N(x,t) = 1 $. 
 	\end{enumerate}
 	Then $ ( X, N, \phi, K, * ) $ is called  fuzzy strong $\phi$-b-normed linear space. 	
  \end{dfn}
  \begin{rem}
If $ K=1 $ and  $ \phi(\alpha) = | \alpha | $ then $ ( X, N,  * ) $ is a B-S type fuzzy normed linear space.
\end{rem}
	%
\begin{eg}\label{ex1}
Consider the linear space $ \mathbb{R} $ and a fuzzy subset $ N $ of $ 	\mathbb{R} \times \mathbb{R} $ by 
\begin{center}
$ N(x,t) = \begin{cases} \frac{t}{t + |x|^p} ~~ t> 0 \\
0~~~~~~~ ~ t \leq 0 
\end{cases} $ 
\end{center}
  for all $  x \in \mathbb{R} $ and $ 0 < p \leq 1 $. \\
Consider the t-norm $ * $ by $ a * b = \min \{ a,b \}, ~ \forall a, b \in \mathbb{R} $. \\
We show that $ N $ is a fuzzy strong $\phi$-b-norm on  $ 	\mathbb{R} \times \mathbb{R} $. For,
\begin{enumerate}[(i)]
\item $ \forall t \in \mathbb{R} $ with $ t \leq 0 $, by definition  we have, $ N(x,t ) = 0 $.   Thus, (bN1) holds.
\item $ ( \forall t \in \mathbb{R}, t > 0, ~ N(x,t) = 1 ) \iff  \frac{t}{t + |x|^p} = 1 
\iff  |x|^p = 0 
\iff  x = 0 $ \\
Therefore (bN2) holds.
\item $ \forall t  > 0 $ and $ c \in \mathbb{R} \setminus \{0 \}  ~$, $ N(cx, t )    = \frac{t}{t + |cx|^p}  = \frac{\frac{t}{|c|^p}}{\frac{t}{|c|^p} + |x|^p}  = N (x, \frac{t}{\phi(c)}) $ \\
where $ \phi(c) = |c|^p, ~ c  \in \mathbb{R}   $ and clearly $ \phi $ satisfies all the conditions of Definition \ref{dfn1}. Thus, (bN3) holds.
\item $ \forall s, t  > 0 $ and $ x,y \in \mathbb{R}  $, $  N(x+y, Ks+ t) = \frac{Ks + t}{Ks +t + |x+y|^p} $ and \\
  $  N(x,s) * N(y,t) = \min \{ N(x,s),  N(y,t) \} = \min \{  \frac{s}{s + |x|^p},  \frac{t}{t + |y|^p} \} $. \\ 
We only prove the inequality for $ s, t > 0 $. \\
Let $  N(x,s) * N(y,t) = \min \{ N(x,s),  N(y,t) \} =  N(x,s) $. \\
Then $ N(y,t) \geq N(x,s) \implies \frac{t}{t + |y|^p} \geq \frac{s}{s + |x|^p} \implies t|x|^p \geq s |y|^p $. \\
Again $ x,y \in \mathbb{R}   $ and $ 0 < p \leq 1 $, 
    $$ |x+y|^p \leq 2^p |x|^p + |y|^p $$
    If we take $ K = 2^p $ then 
    \begin{align*}
    N(x+y, 2^ps+ t) -  N(x,s) & =  \frac{2^ps + t}{2^ps +t + |x+y|^p} - \frac{s}{s + |x|^p} \\
    & \geq   \frac{2^ps + t}{2^ps +t + { 2^p |x|^p + |y|^p}} - \frac{s}{s + |x|^p} \\ 
    & = \frac{ t|x|^p - s |y|^p} {(2^ps +t + { 2^p |x|^p + |y|^p}) (s + |x|^p)} \geq 0
    \end{align*}
    Hence $ N(x+y, 2^ps+ t) \geq   N(x,s) = \min \{ N(x,s),  N(y,t) \} $. \\
    Similarly,  it can be shown that if $  \min \{ N(x,s),  N(y,t) \} =  N(y,t) $ then \\
    $ N(x+y, 2^ps+ t) \geq   N(y,t) = \min \{ N(x,s),  N(y,t) \} $. \\
    Therefore, (bN4) holds.
\item From the definition of $ N(x,t) $ it is clear that $ N(x,.) $ is a non-decreasing function of $ t $ and 
 $ \underset{t \rightarrow \infty }{lim} N(x,t) = 1 $. 
\end{enumerate}

Hence $ ( X, N, \phi, K, * ) $ is a  fuzzy strong $\phi$-b-normed linear space where $ K = 2^p(>1) $ and \\ 
$ \phi(\alpha) = |\alpha|^p, \forall \alpha
  \in \mathbb{R}, ~ 0< p\leq 1 $. 
\end{eg}
\begin{eg}\label{ex2}
Consider the linear space $ \mathbb{R} $ and a fuzzy subset $ N $ of $ 	\mathbb{R} \times \mathbb{R} $ by 
\begin{center}
$ N(x,t) = \begin{cases} \exp({-\frac{|x|^p}{t }}) ~~~ t> 0 \\
0 ~~~~~~~~~~~~~~~ t \leq 0 
\end{cases} $ 
\end{center}
for all$  x \in \mathbb{R} $ and $ 0 < p \leq 1 $ and consider the t-norm $ * $ by $ a * b =  ab , ~ \forall a, b \in \mathbb{R} $. Now,
\begin{enumerate}[(i)]
\item Clearly (bN1) holds from the definition.
\item  $ ( \forall t \in \mathbb{R}, t > 0, ~ N(x,t) = 1 ) \iff  \exp({-\frac{|x|^p}{t }}) = 1 \iff  |x|^p = 0 \iff  x = 0 $ \\
Therefore (bN2) holds.
\item $ \forall t  > 0 $ and $ c \in \mathbb{R}\setminus \{0\}  $, 
$~~ N(cx, t )    =  \exp({-\frac{|cx|^p}{t }})    =    \exp({-\frac{|x|^p}{\frac{t}{|c|^p} }})        = N (x, \frac{t}{\phi(c)}) $ \\
where $ \phi(c) = |c|^p, ~ c  \in \mathbb{R}   $ and clearly $ \phi $ satisfies all the conditions of Definition \ref{dfn1}. Thus, (bN3) holds.
\item For $ s, t  > 0 $ and $ x,y \in \mathbb{R}  $, $  N(x+y, Ks+ t) = \exp({-\frac{|x + y |^p}{ Ks + t }})  $ and \\
  $  N(x,s) * N(y,t) =  N(x,s) \cdot  N(y,t)  =   \exp({-\frac{|x|^p}{s }}) \cdot \exp({-\frac{|y|^p}{t }}) $. \\ 
Using the inequality, $|x+y|^p \leq 2^p |x|^p + |y|^p $,    $ ~ x,y \in \mathbb{R} $ and $ 0 < p \leq 1 $ and taking  $ K = 2^p $, we obtain 
$$ - \frac{|x + y |^p}{ 2^ps + t }  \geq - \frac{ 2^p |x|^p + |y|^p }{ 2^ps + t }  \geq - \frac{ 2^p |x|^p }{ 2^ps + t } ~ - \frac{ |y|^p }{ 2^ps + t }  \geq - \frac{ |x|^p }{ s  } ~ - \frac{ |y|^p }{ t } $$
 which implies $  N(x+y, 2^ps+ t) \geq   N(x,s) \cdot  N(y,t) $. \\
 Thus (bN4): $  N(x+y, 2^ps+ t) \geq   N(x,s) *  N(y,t) $ holds $ \forall s, t \in \mathbb{R} $ and $\forall x,y \in \mathbb{R}  $.
 \item Clearly $ N(x,\cdot) $  is a non-decreasing function of $ t $ and 
 $ \underset{t \rightarrow \infty }{lim} N(x,t) = 1 $. 
\end{enumerate}
Hence $ ( X, N, \phi, K, * ) $ is a  fuzzy strong $\phi$-b-normed linear space where $ K = 2^p(>1) $ and \\ 
$ \phi(\alpha) = |\alpha|^p, \forall \alpha
  \in \mathbb{R}, ~ 0< p\leq 1 $.	 
\end{eg} 
\begin{rem}
The notions of converges, Cauchy sequences, boundedness, etc. are same as definitions in Bag and Samanta type fuzzy normed linear space\cite{[1]}.
\end{rem}
\section{Finite dimensional fuzzy strong $\phi$-b-normed linear spaces}
 In this section some basic results on finite dimensional fuzzy strong $\phi$-b-normed linear spaces are established.
 \begin{lma}\label{lma1}
 Let $ ( X, N, \phi, K, * ) $ be a  fuzzy strong $\phi$-b-normed linear space with the underlying $t$-norm $*$ continuous at $ (1,1) $ and 
 $ \{ x_1, x_2, \cdots x_n \} $ be a linearly independent set of vectors in $ X $. Then $ \exists c > 0 $ and $ \delta \in ( 0, 1 ) $ such that 
 for any set of scalars  $ \{ \alpha_1, \alpha_2, \cdots \alpha_n \} $ with $ \sum _{i=1} ^n |\alpha_i| \neq 0 $,  
 
 \begin{equation}\label{eqn1}
   N( \alpha_1x_1 +  \alpha_2x_2 + \cdots + \alpha_nx_n, \frac{Kc}{\phi( \frac{1}{\sum _{i=1} ^n |\alpha_i|} ) } ) < 1 - \delta 
 \end{equation}
 \begin{proof}
 The relation (\ref{eqn1}) is equivalent to the relation 
 \begin{equation}\label{eqn2}
  N( \beta_1x_1 +  \beta_2x_2 + \cdots + \beta_nx_n, Kc ) < 1 - \delta 
\end{equation}
for some $ c > 0 $ and $ \delta \in ( 0, 1 ) $ and for all set of scalars   $ \{ \beta_1, \beta_2, \cdots, \beta_n \} $ with 
$ \sum _{i=1} ^n |\beta_i| = 1 $. \\
If possible suppose that (\ref{eqn2}) does not hold. Thus for each $ c > 0 $ and $ \delta \in ( 0, 1 ), \exists $ a set of scalars 
  $ \{ \beta_1, \beta_2, \cdots, \beta_n \} $ with $ \sum _{i=1} ^n |\beta_i| = 1 $   for which 
  \begin{equation*}
  N( \beta_1x_1 +  \beta_2x_2 + \cdots + \beta_nx_n, Kc ) \geq 1 - \delta 
\end{equation*}
Then for $ c = \delta = \frac{1}{m}, ~ m =1, 2, \cdots, ~\exists $ a set of scalars 
  $ \{ \beta_1^{(m)}, \beta_2^{(m)}, \cdots, \beta_n^{(m)} \} $ with $ \sum _{i=1} ^n |\beta_i^{(m)}| = 1 $ such that  
  $ N(y_m, \frac{K}{m} )  \geq 1 - \frac{1}{m} $ where $ y_m = \beta_1^{(m)}x_1 +  \beta_2^{(m)}x_2 + \cdots + \beta_n^{(m)}x_n $. \\
  Since $ \sum _{i=1} ^n |\beta_i^{(m)}| = 1 $, we have $ 0 \leq |\beta_i^{(m)}| \leq 1 $ for $ i = 1,2, \cdots , n $. So for each fixed $ i $, 
  the sequence $ \{ \beta_i^{(m)} \} $ is bounded and hence $ \{ \beta_i^{(m)} \} $  has a convergent subsequence. Let $ \beta_1 $ denotes the
  limit of that subsequence and let $ \{ y _{1,m} \} $ denotes the corresponding subsequence of $ \{ y_m \} $. \\
  By the same argument  $ \{ y _{1,m} \} $ has a subsequence $ \{ y _{2,m} \} $ for which the corresponding subsequence of scalars 
  $ \{ \beta_2^{(m)} \} $ converges to $ \beta_2 $.\\
   Continuing in this way, after $ n $ steps we obtain a subsequence $ \{ y _{n,m} \} $ where $ y _{n,m} = \sum _{i=1} ^n \gamma_i ^{(m)} x_i  $ 
   with  $ \sum _{i=1} ^n |\gamma_i^{(m)}| = 1 $ and $ \gamma_i^{(m)} \rightarrow \beta _i $ as $ m \rightarrow \infty $ for each 
   $ i = 1, 2, \cdots, n $.\\
   Let $ y = \beta_1x_1 +  \beta_2x_2 + \cdots + \beta_nx_n $. Now, 
   \begin{align*}
   N(y_{n,m} - y, t) & = N(\sum_{j=1}^n(\gamma_j^{(m)} - \beta_j)x_j, t )\\
   & =  N( (\gamma_1^{(m)} - \beta_1)x_1 + \sum_{j=2}^n(\gamma_j^{(m)} - \beta_j)x_j, \frac{t}{n}+ K(n-1) \frac{t}{nK} )\\
   & \geq N( (\gamma_1^{(m)} - \beta_1)x_1, \frac{t}{n} ) * N ( \sum_{j=2}^n(\gamma_j^{(m)} - \beta_j)x_j, (n-1) \frac{t}{nK} ) \\
   & = N( (\gamma_1^{(m)} - \beta_1)x_1, \frac{t}{n} ) * N ( (\gamma_2^{(m)} - \beta_2)x_2 + \sum_{j=3}^n(\gamma_j^{(m)} - \beta_j)x_j, 
   \frac{t}{nK} + K(1-\frac{2}{n}) \frac{t}{K^2} ) \\
    & \geq N( (\gamma_1^{(m)} - \beta_1)x_1, \frac{t}{n} ) * N ( (\gamma_2^{(m)} - \beta_2)x_2, \frac{t}{nK} ) *
     N (   \sum_{j=3}^n(\gamma_j^{(m)} - \beta_j)x_j,  (1-\frac{2}{n}) \frac{t}{K^2} ) \\
      & \cdots \\
      & \geq N( (\gamma_1^{(m)} - \beta_1)x_1, \frac{t}{n} ) * N ( (\gamma_2^{(m)} - \beta_2)x_2, \frac{t}{nK} ) * \cdots * 
      N ( (\gamma_n^{(m)} - \beta_n)x_n, \frac{t}{nK^{n-1}} ) \\
    & =  N( x_1, \frac{t}{n \phi((\gamma_1^{(m)} - \beta_1))} )  * \cdots * 
      N ( x_n, \frac{t}{nK^{n-1} \phi((\gamma_n^{(m)} - \beta_n))} ) \\     
   \end{align*}
   Taking limit as $ m \rightarrow \infty $ on both sides, we have $ \underset{ m \rightarrow \infty }{lim} N(y_{n,m} - y, t)  \geq 1 * 1* \cdots * 1, 
   ~ \forall t > 0 $ which implies $  \underset{ m \rightarrow \infty }{lim} N(y_{n,m} - y, t)  = 1,  ~ \forall t > 0 $.\\
   Now for $ r > 0 $, choose $ m $ such that $ \frac{1}{m} < \frac{r}{K^2} $. We have  
    \begin{align*}
  &   N(y_{n,m}, \frac{r}{K} )  = N ( y_{n,m} + \theta, \frac{K}{m} + K(\frac{r}{K^2} - \frac{1}{m}) )  \geq  N ( y_{n,m}, \frac{K}{m} ) *
   N(  \theta,   (r-\frac{1}{m}) ) 
      \geq  ( 1 - \frac{K}{m} ) * 1\\
   \implies  & \underset{ m \rightarrow \infty }{lim} N(y_{n,m}, \frac{r}{K})  \geq 1  \\
   \implies &  \underset{ m \rightarrow \infty }{lim} N(y_{n,m}, \frac{r}{K} )  = 1  
   \end{align*}
   Again, 
   \begin{align*}
    & N(y, 2r )  = N ( y- y_{n,m} + y_{n,m} , r + K \cdot \frac{r}{K}  )  \geq  N ( y- y_{n,m}, r ) * N(  y_{n,m}, \frac{r}{K}    )  \\
     \implies & N(y, 2r ) \geq  \underset{ m \rightarrow \infty }{lim}   N ( y- y_{n,m}, r ) *  \underset{ m \rightarrow \infty }{lim} 
      N(  y_{n,m}, \frac{r}{K}    )   \\
      \implies & N(y, 2r ) \geq 1 * 1 =   1 \\
      \implies  & N(y, 2r ) = 1
   \end{align*}  
   Since $ r > 0 $ is arbitrary, so $ y =\theta $. \\
   Again since $ \sum _{i=1} ^n |\beta_i^{(m)}| = 1 $ and $ \{ x_1, x_2, \cdots x_n \} $ is a linearly independent set of vectors so 
   $ y = \beta_1x_1 +  \beta_2x_2 + \cdots + \beta_nx_n  \neq \theta $. \\
   Thus we arrive at a contradiction. Hence (\ref{eqn2}) holds and Lemma is proved.
 \end{proof}
 \end{lma}
 \begin{thm}\label{thm1}
 Every finite dimensional fuzzy strong $\phi$-b-normed linear space with the underlying $t$-norm $*$ continuous at $ (1,1) $ is complete.
 \begin{proof}
 Let $ ( X, N, \phi, K, * ) $ be a  fuzzy strong $\phi$-b-normed linear space  where $ K (> 1) $ is a real constant. \\
 Let $ dimX= r $ and $ \{ e_1, e_2, \cdots, e_r \} $ be a basis for $ X $. \\
 Let $ \{ x_p \}  $ be a Cauchy sequence in $ X $. Then $ x_n = \beta_1^{(n)}e_1 +  \beta_2^{(n)}e_2 + \cdots +  \beta_r^{(n)}e_r $ for suitable 
 scalars $  \beta_1^{(n)},  \beta_2^{(n)}, \cdots, \beta_r^{(n)} $. So,
 \begin{equation}\label{eqn3}
 \underset{ m, n \rightarrow \infty }{lim}   N ( x_m - x_n, t ) = 1, ~ \forall t > 0
\end{equation}  
 Now from Lemma \ref{lma1}, it follows that    $ \exists c > 0 $ and $ \delta \in ( 0, 1 ) $ such that
 \begin{equation}\label{eqn4}
 N ( \sum_{i=1}^r {(\beta_i^{(m)} - \beta_i^{(n)}) e_i}, \frac{cK}{ \phi ( \frac{1}{ \sum_{i=1}^r {|\beta_i^{(m)} - \beta_i^{(n)}| }} )} ) < 1 - \delta
 \end{equation}
 If $ \sum_{i=1}^r {|\beta_i^{(m)} - \beta_i^{(n)}|} = 0 $ then $ \beta_i^{(m)} = \beta_i^{(n)}, ~ \forall i $ implies $ \{x_n\} $ is a constant sequence and hence follows the theorem. So we may assume $ \sum_{i=1}^r {|\beta_i^{(m)} - \beta_i^{(n)}|} \neq 0 $. \\
 Again for $ 0 < \delta < 1 $, from (\ref{eqn3}), it follows that there exist a positive integer $ n_0(\delta, t) $ such that 
  \begin{equation}\label{eqn5}
 N ( \sum_{i=1}^r {(\beta_i^{(m)} - \beta_i^{(n)}) e_i}, t) > 1 - \delta, ~ \forall m, n \geq n_0(\delta, t)
 \end{equation}
 Now from (\ref{eqn4}) and (\ref{eqn5}) we have, 
 \begin{align*}
 &N ( \sum_{i=1}^r {(\beta_i^{(m)} - \beta_i^{(n)}) e_i}, \frac{cK}{ \phi ( \frac{1}{ \sum_{i=1}^r {|\beta_i^{(m)} - \beta_i^{(n)}| }} )} ) < 
 N ( \sum_{i=1}^r {(\beta_i^{(m)} - \beta_i^{(n)}) e_i}, t), ~ \forall m, n \geq n_0(\delta, t)   \\
 \implies & \frac{cK}{ \phi( \frac{1}{  \sum_{i=1}^r {|\beta_i^{(m)} - \beta_i^{(n)}| })} } < t, ~ \forall m, n \geq n_0(\delta, t) ~~ 
( \text{Since} ~ N(x,\cdot) ~ \text {is non-decreasing})  \\
 \end{align*}
 Since $ t > 0$ is arbitrary, thus 
 \begin{align*}
 & \underset{ m, n \rightarrow \infty }{lim} \frac{cK}{ \phi ( \frac{1}{ \sum_{i=1}^r {|\beta_i^{(m)} - \beta_i^{(n)}| }} )} = 0 \\
 \implies & \underset{ m, n \rightarrow \infty }{lim} \phi ( \frac{1}{ \sum_{i=1}^r {|\beta_i^{(m)} - \beta_i^{(n)}| }} ) = \infty \\
 \implies &  \phi ( \frac{1}{ \underset{ m, n \rightarrow \infty }{lim}  \sum_{i=1}^r {|\beta_i^{(m)} - \beta_i^{(n)}| }} ) = \infty ~~ 
( \text{Since} ~ \phi ~ \text {is continuous }) \\
\implies &  \underset{ m, n \rightarrow \infty }{lim}  \sum_{i=1}^r {|\beta_i^{(m)} - \beta_i^{(n)}| } = 0 ~~ ( \text{Since} ~
\underset{ \alpha \rightarrow \infty }{lim} \phi(\alpha) = \infty ) ~
 \end{align*}
 Therefore, $ \{ \beta_i^{(m)} \} $ is a Cauchy sequence of scalars for each $ i = 1, 2, \cdots, r $. So each sequence $ \{ \beta_i^{(m)} \} $
  converges. \\
  Let $ \underset{  n \rightarrow \infty }{lim} \beta_i^{(n)} = \beta_i  $ for $ i = 1, 2, \cdots, r $.  Define $ x = \sum_{i=1}^r \beta_ie_i $.
    Clearly $ x \in X $. \\
    By similar calculation as in Lemma \ref{lma1}, it can be shown that $  \underset{  n \rightarrow \infty }{lim} N (x_n -x,t) = 1, ~ \forall t > 0 $.\\
    Hence $X$ is complete.
 \end{proof}
 \end{thm}
 \begin{thm}\label{thm2}
 Let $ ( X, N, \phi, K, * ) $ be a finite dimensional  fuzzy strong $\phi$-b-normed linear space in which the underlying $t$-norm $*$ continuous 
 at $ (1,1) $. Then a subset $ A $ of $X $ is compact iff $ A $ is closed and bounded. 
 \begin{proof}
 First we suppose that $ A $ is compact. We have to show that $ A $ is closed and bounded. \\
 For, let $ x \in \bar{A} $. Then there exist a sequence $ \{ x_n \} $ in $ A $ such that $ \underset{ n \rightarrow n }{lim} x_n = x $. \\
 Since $ A $ is compact, there exist a subsequence $ \{ x_{n_k} \} $ of  $ \{ x_n \} $ converges to a point in $ A $. Again $ x_n \rightarrow x $ 
 so $ x_{n_k} \rightarrow x $ and hence $ x \in A $. So $ A $ is closed. \\
 If possible suppose that $ A $ is not bounded. \\
 Then $ \exists r_0, ~ 0<r_0 < 1 $, such that for each positive integer $ n, ~ \exists x_n \in A $ for which $ N (x_n,n ) \leq 1 - r_0 $. \\
 Since $ A $ is compact, there exist  a subsequence $ \{ x_{n_p} \} $ of  $ \{ x_n \} $ converging to some element $ x \in A $. Thus 
 $ \underset{ p \rightarrow \infty}{lim} N(x_{n_p} - x , t ) = 1, ~ \forall t>0 $. \\
 Again,  $ N (x_{n_p},n_p ) \leq 1 - r_0 $. Now,
 $$ 1 - r_0 \geq N (x_{n_p},n_p ) =  N (x_{n_p} - x + x, \frac{t}{K} + K (\frac{n_p}{K} - \frac{t}{K^2} )) 
 \geq  N (x_{n_p} - x, \frac{t}{K} ) * N ( x,  (\frac{n_p}{K} - \frac{t}{K^2} )) $$
 On above inequality, taking limit as $ p \rightarrow \infty $, we obtain
 \begin{align*}
&  1 - r_0 \geq \underset{ p \rightarrow \infty}{lim} N (x_{n_p} - x, \frac{t}{K} ) * \underset{ p \rightarrow \infty}{lim} 
 N ( x,  (\frac{n_p}{K} - \frac{t}{K^2} )) \\
 \implies & 1 - r_0 \geq 1* 1 = 1 \\
 \implies & r_0 \leq 0
 \end{align*} 
 This is a contradiction. Hence $A$ is bounded. \\
 Conversely suppose that $ A $ is closed and bounded and we have to show that $ A $ is compact. \\
 Let $ dimX= r $ and $ \{ e_1, e_2, \cdots, e_r \} $ be a basis for $ X $. \\
 Let us choose a sequence $ \{ x_p \} $ in $ A $ and suppose $ x_p = \beta_1^{(p)}e_1 +  \beta_2^{(p)}e_2 + \cdots +  \beta_r^{(p)}e_r $ 
 for suitable  scalars $  \beta_1^{(p)},  \beta_2^{(p)}, \cdots, \beta_r^{(p)} $.  \\
 Now from Lemma \ref{lma1},  $ \exists c > 0 $ and $ \delta \in ( 0, 1 ) $ such that
 \begin{equation}\label{eqn6}
   N( \sum_{i=1}^r \beta_i^{(p)}e_i, \frac{Kc}{\phi( \frac{1}{\sum _{i=1} ^n |\beta_i^{(p)}|} ) } ) < 1 - \delta 
 \end{equation}
 Again since $A $ is bounded, for $ \delta \in ( 0,1), ~ \exists t > 0 $ such that $ N (x, t ) > 1 - \delta, \forall x \in A $. So
 \begin{equation}\label{eqn7}
   N( \sum_{i=1}^r \beta_i^{(p)}e_i, t ) > 1 - \delta 
 \end{equation}
 From (\ref{eqn6}) and (\ref{eqn7}) we get, 
 \begin{align*}
 & N (\sum_{i=1}^r \beta_i^{(p)}e_i, \frac{Kc}{\phi( \frac{1}{\sum _{i=1} ^n |\beta_i^{(p)}|} ) }  ) < 
 N ( \sum_{i=1}^r \beta_i^{(p)}e_i, t )  \\
 \implies & \frac{Kc}{\phi( \frac{1}{\sum _{i=1} ^n |\beta_i^{(p)}|} ) } < t, ~ \forall m, n \geq n_0(\delta, t) ~~ 
( \text{Since} ~ N(x,\cdot) ~ \text {is non-decreasing})  \\
 \end{align*}
 Without loss of generality we may assume that $ \sum _{i=1} ^n |\beta_i^{(p)}| \neq 0 $. \\
 If $ \sum _{i=1} ^n |\beta_i^{(p)}| = 0 $ then $ \beta_i^{(p)} = 0 $, for $ i = 1, 2, \cdots $ and $ \forall p $. Then $ \{ x_p \} $ is a 
 constant sequence and the theorem follows. \\
 Since $ c, K, t $ are three fixed positive real numbers, it follows that $ 0 <  \sum _{i=1} ^n |\beta_i^{(p)}|  < \infty $. \\
 Therefore the sequence of scalars $ \beta_i^{(p)}, p = 1, 2, \cdots  $ and for $ i = 1,2,\cdots, n $ is bounded. So  by Bolzano-Weierstrass theorem,
  there exist a convergent subsequence of $ \{ \beta_i^{(p)}\} $. Now, we follow the techniques of Lemma \ref{lma1} to show that there exist 
  a subsequence of  $ \{ x_p \} $ that converges to some point in $ A $. \\
  Thus $ A $ is compact and this proves the theorem.
 \end{proof}
 \end{thm}
 %
 
%



~\\
~\\
%
%
    \textbf{Conclusion:} Recently different types of generalized fuzzy metric spaces as well as  generalized fuzzy normed linear spaces have been
     developed by several authors. Following the definition of fuzzy strong b-metric spaces, we introduce the idea of   fuzzy strong $\phi$-b-normed 
     linear spaces and study some results in finite finite dimensional fuzzy strong $\phi$-b-normed  linear spaces. We think there is a huge scope of
      research to develop fuzzy strong $\phi$-b-normed  linear spaces. Results on completeness and compactness, operator norms etc. are the open 
      problems in such spaces.   \\
~\\
    \textbf{Acknowledgment:}
 The author AD is thankful to University Grant Commission (UGC), New Delhi, India for awarding her senior research fellowship [Grant No.1221/(CSIRNETJUNE2019)]. We are also grateful to Department of Mathematics, Siksha-Bhavana, Visva-Bharati. \\

\end{document}